\documentclass[psamsfonts]{amsart}

\newtheorem{theorem}{Theorem}[]
\newtheorem{lemma}[theorem]{Lemma}

\newtheorem{corollary}[theorem]{Corollary}

\theoremstyle{definition}

\newtheorem{example}[theorem]{Example}

\newtheorem{remark}[theorem]{Remark}
\newtheorem{remarks}[theorem]{Remarks}

\theoremstyle{remark}

\numberwithin{equation}{section}

\newcommand{\ZZ}{\mathbb{Z}}

\newcommand  {\foa}     {\mathfrak{a}}

\newcommand  {\fop}     {\mathfrak{p}}

\newcommand  {\foq}     {\mathfrak{q}}

\newcommand  {\fos}     {\mathfrak{r}}

\newcommand  {\ra}      {\rightarrow}

\newcommand  {\Spec}    {\operatorname{Spec}}

\def\mydate{\number\day\space\ifcase\month \or January\or February\or March\or April\or May\or
June\or July\or August\or September\or October\or November\or
December\fi \space\number\year}

\usepackage{amscd}
\usepackage{amssymb}

\begin{document}

\title[Lifting chains of prime ideals]
{Lifting chains of prime ideals}


\author[Holger Brenner]{Holger Brenner}
\address{Mathematische Fakult\"at, Ruhr-Universit\"at, Universit\"atsstr. 150,
               44780 Bochum, Germany}
\email{brenner@cobra.ruhr-uni-bochum.de}





\begin{abstract}
We give an elementary proof that for a ring homomorphism $A \ra B$
satisfying the property that every ideal in $A$ is contracted from $B$
the following property holds:
for every chain of prime ideals
$ \fop_0 \subset  \ldots  \subset \fop_r $ in $A$ there exists a
chain of prime ideals
$\foq_0 \subset  \ldots  \subset \foq_r $ in $B$ such that
$\foq_i \cap A= \fop_i$.
\end{abstract}

\maketitle


\bigskip
Let $A$ and $B$ be commutative rings and let
$\varphi: A \ra B$ be a ring homomorphism.
This induces a continouus mapping
$\varphi^*: \Spec \, B \ra \Spec\, A$ by sending a prime ideal
$\foq \subset B$ to $\varphi^{-1}(\foq)$.
Properties of the ring homomorphism are then often reflected
by topological properties of $\varphi^*$.
For example, if $A \ra B$ is integral, then ``going up''
holds,
and if $A \ra B$ is flat, then ``going down'' holds
(see \cite[Proposition 4.15 and Lemma 10.11]{eisenbud}.
If moreover $\varphi ^*:\, \Spec \, B \ra \Spec \,A$
is surjective and going up or going down holds,
then also the following property holds:
for every given chain of prime ideals
$ \fop_0 \subset  \ldots  \subset \fop_r $ in $A$ there exists a
chain of prime ideals
$\foq_0 \subset  \ldots  \subset \foq_r $ in $B$ lying over it.

In this note we give a direct and elementary proof
showing that this chain lifting property holds also under the condition
that every ideal in $A$ is contracted from $B$, i.e.
$I= \varphi^{-1}(IB)$ holds for every ideal $I \subseteq A$.
This result can be found for pure homomorphisms in
Picavet's paper (see \cite{picavet}[Proposition 60 and Theorem 37])
and is proved using valuation theory.
Our direct method allows to find explicitely chains of prime ideals and
characterizes which prime ideals $\foq_0$ over $\fop_0$ may be extended
to a chain.
We start with the following lemma.

\begin{lemma}
\label{lemma}
Let $B$ be a commutative ring, let
$\foa_0 ,  \ldots ,\foa_r $ be ideals and
$F_0, \ldots ,F_r$ multiplicatively closed systems.
Define inductively {\rm (}set $S_{r+1}=\{1\}${\rm )}
for $i=r, \ldots, 0$ the following multiplicatively closed sets
$$ S_i 
= \{ s \in B:\, (s, { \foa }_i) \cap F_i \cdot S_{i+1} \neq \emptyset \} \, .$$
Then the following are equivalent.

\renewcommand{\labelenumi}{(\roman{enumi})}
\begin{enumerate}

\item
$0 \not\in S_0$.

\item
$ { \foa }_i \cap F_i \cdot S_{i+1} = \emptyset $
for
$i= 0, \ldots , r$.

\item
There exists a
chain of prime ideals
${ \foq}_0 \subseteq  \ldots  \subseteq { \foq}_r $
such that
${ \foa }_i \subseteq { \foq}_i $
and
${ \foq}_i \cap F_i \cdot S_{i+1} = \emptyset $.

\item
There exists a
chain of prime ideals
${ \foq}_0 \subseteq  \ldots  \subseteq { \foq}_r $
such that
${ \foa }_i \subseteq { \foq}_i $
and
${ \foq}_i \cap F_i = \emptyset $.
\end{enumerate}
\end{lemma}

\proof
It is clear that the $S_i$ are multiplicatively closed
and that $S_{i+1} \subseteq S_i$.
(i) $\Leftrightarrow $ (ii).
If $0 \in S_0$, then
$\foa_0 \cap F_0 \cdot S_{i+1} \neq \emptyset$, and if
$\foa_i \cap F_i \cdot S_{i+1} \neq \emptyset$ for some $i$,
then $0 \in S_i$ and thus also $0 \in S_0$.

We show (ii) $ \Rightarrow $ (iii) by induction.
Since
$ \foa_0 \cap F_0 S_1 = \emptyset $,
there exists (\cite[Ch.2 \S5, Corollary 2]{bourbakiac}) a prime ideal $\foq_0$
such that
$\foa_0 \subseteq \foq_0$ and $ \foq_0 \cap F_0S_1 = \emptyset $.

Thus suppose that the chain
$\foq_0 \subset \ldots \subset \foq_i$ is already constructed.
We have to look for a prime ideal $\foq_{i+1}$ which includes
both ${ \foq}_i$ and ${ \foa }_{i+1}$ and which is disjoint
to $F_{i+1} \cdot S_{i+2} $.
If such a prime ideal would not exist,
then
$({ \foq}_i + { \foa }_{i+1}) \cap F_{i+1} \cdot S_{i+2} \neq \emptyset $,
say $q+a=f \cdot s $, where $q \in \foq_i$, $a \in \foa_{i+1}$,
$f \in F_{i+1}$ and $s \in S_{i+2}$.
Then by definition $q \in S_{i+1} $ contradicting the induction
assumption.

(iii) $\Rightarrow $ (iv) and (iii) $\Rightarrow $ (ii) are clear,
so we have to show (iv) $\Rightarrow$ (iii).
We show this by descending induction, the beginning for $i=r$ is clear.
Suppose that $\foq_{i-1} \cap F_{i-1}S_i \neq \emptyset$,
and let $q=fs$ be an element in the intersection,
$q \in \foq_{i-1}$, $f\in F_{i-1}$, $s \in S_i$.
Since $F_{i-1}$ is disjoined to the prime ideal $\foq_{i-1}$, it follows
that
$s \in \foq_{i-1}$.
On the other hand, since $s \in S_i$ we have an equation
$bs +q' =f's'$, where $b \in B$,
$q' \in \foq_{i}$, $f' \in F_{i}$, $s' \in S_{i+1}$,
and this  contradicts the induction hypothesis.
\qed

\begin{remark}
The referee (whom I thank for his careful reading)
pointed out that there exists a similar and more general
result in a preprint of G. Bergman (see \cite{bergman}).
Bergman studies for a partially ordered set $I$
and ideals $\foa_i$ and multiplicatively closed subsets $S_i$
in a commutative ring
the existence of prime ideals $\fop_i$, $\foa_i \subseteq \fop_i$,
$\fop_i \cap S_i = \emptyset$ such that $\fop_i \subset \fop_j$
holds for $i \leq j$. \cite[Proposition 9]{bergman}
gives a characterization for the existence of such prime ideals
for a tree order $I$ in terms of an inductively
defined system of equations which is related to our
characterization in Lemma \ref{lemma}(ii). It is possible that using Bergman's
result one may obtain a stronger version of the following
theorem.
\end{remark}

\begin{theorem}
Let $A$ and $B$
be commutative rings and let $\varphi:A \ra B$ be a
ring homomorphism such that $I=\varphi^{-1} (IB)$ holds
for every ideal $I \subseteq A$.
Then for every chain of prime ideals
${ \fop}_0 \subset  \ldots  \subset { \fop}_r $ in $\Spec\, A$
there exists a chain of prime ideals
${ \foq}_0 \subset  \ldots  \subset { \foq}_r $
in $B$ such that $\fop_i= \foq_i \cap A$ for $i=0, \ldots, r$.
\end{theorem}

\proof
Let a chain of prime ideals
${ \fop}_0 \subset  \ldots  \subset { \fop}_r $
in $A$ be given.
We shall apply the preceeding lemma
to the ideals $\foa_i=  \fop_i B$
and the multiplicatively closed sets
$F_i = A - { \fop}_i \subset B$.
Note that the fiber over $\fop$ consists of the prime ideals $\foq$
for which $\fop B \subset \foq$ and $\foq \cap \varphi(A- \fop) = \emptyset$
hold.
Define $S_i \subseteq B$ as before and suppose that $0 \in S_0$.
This means that
there exists an element $a_0 \in \foa_0$ such that
$a_0 = f_0 \cdot s_1 $, where $f_0 \in F_0$, $s_1 \in S_1$.
This means by definition that we have an equation
$$ b_1s_1 + a_1 =f_1s_2 ,\,
\mbox{ where } \, b_1 \in B, \, a_1 \in { \foa }_1, \,
f_1 \in F_1 \, \mbox{ and }\, s_2 \in S_2 \, .$$
Going on recursively we find equations
$$ b_j s_j+a_j=f_js_{j+1}, \, \mbox{ where }\,
b_j \in B,\, a_j \in { \foa }_j,\, f_j \in F_j \, \mbox{ and }\,
s_{j+1} \in S_{j+1} \, , $$
and eventually
$$b_r s_r + a_r = f_r, \, \mbox{ where } \,
b_r,s_r \in S_r,\, a_r \in { \foa}_r ,\, f_r \in F_r \, .$$
We multiply the last equation by
$f_{r-1} \cdots f_0 $ and get
$$b_r(s_rf_{r-1})f_{r-2} \cdots f_0 +a_r f_{r-1} \cdots f_0
=f_r f_{r-1} \cdots  f_0  \, .$$
We may replace $b_r(s_rf_{r-1})f_{r-2} \cdots f_0$ by
$$b_r(b_{r-1}s_{r-1} +a_{r-1})f_{r-2} \cdots f_0 =
b_rb_{r-1}(s_{r-1}f_{r-2}) \cdots f_0
+b_ra_{r-1}f_{r-2} \cdots f_0 \, ,$$
and so going on we find that
$f_r \cdots f_0 =$
$$
b_r \cdots b_1a_0 + b_r \cdots b_2a_1f_0 +b_r \cdots b_3a_2f_1f_0
+ \ldots +
b_ra_{r-1}f_{r-2} \cdots f_0
+a_rf_{r-1} \cdots f_0 \,.$$
This equation shows that
$$ f_r \cdots f_0 \in
(\fop_0 + \fop_1 f_0 + \fop_2f_1f_0
+ \ldots +
\fop_{r-1}f_{r-2} \cdots f_0
+ \fop_rf_{r-1} \cdots f_0)B \, $$
and this yields an equation in $A$ (here we apply the condition
that every ideal is contracted),
$$p_0+ p_1 f_0+p_2 f_1f_0+ \ldots +p_{r-1}f_{r-2} \cdots f_0+
p_rf_{r-1} \cdots f_0 -f_r \cdots f_0 =0 \, ,$$
where $p_i \in \fop_i$.
We may write this as
$$p_0 =-f_0
( p_1+p_2f_1+ \ldots +p_{r-1}f_{r-2} \cdots f_1+
p_rf_{r-1} \cdots f_1 -f_r \cdots f_1)$$
and therefore
$p_1+p_2f_1+ \ldots +p_{r-1}f_{r-2} \cdots f_1+
p_rf_{r-1} \cdots f_1 -f_r \cdots f_1 \in \fop_0 \subset \fop_1$.
Then again we may multiply out $f_1$ and so on until we find
$p_{r-1}+ p_rf_{r-1}-f_rf_{r-1} \in \fop_{r-2} \subset \fop_{r-1}$
and then $p_rf_{r-1} -f_rf_{r-1} \in \fop_{r-1}$,
hence $p_r-f_r \in \fop_{r-1}$ and $f_r \in \fop_r$,
which is a contradiction.
\qed

\begin{remarks}
The condition that every ideal is contracted is fulfilled for example if
$ \varphi:A \ra B$ is a pure homomorphism. This
means that for every $A$-module
$M$ the natural mapping $M \ra M \otimes_A B$ is injective.
If $B$ contains $A$ as a direct summand,
then $A \subseteq B$ is pure.
Direct summands arise often in invariant theory:
if a linearly reductive group acts on a ring $B$,
then the ring of invariants
$A=B^G$ is a direct summand in $B$, see \cite[Ch.1, \S1]{GIT}.
Example \ref{example} below shows that for a direct summand
neither going up nor going down hold in general.

G. Picavet studies in \cite{picavet} the property of a ring homomorphism
that over every chain of prime ideals $\fop \subset \foq $
there lies a chain of prime ideals above. He calls a ring homomorphism with
this property subtrusif and shows that a homomorphism $\varphi :A \ra B$
is universally
subtrusif if and only if for every valuation domain $A \ra V$ the
corresponding homomorphism $V \ra B \otimes_A V$ is pure.

Picavet proved the theorem for universally subtrusive
morphisms
\cite[Proposition 60 in connection with Theorem 37]{picavet}
using several facts from valuation theory:
that for a chain of prime ideals $\fop_0 \subset \ldots  \subset \fop_r$ in
a domain $A$ there exists a valuation ring $A \subseteq V \subseteq Q(A)$
and a chain of prime ideals $\fos_0 \subset \ldots \subset \fos_r$ in $V$
with $\fos_i \cap A = \fop_i$, see \cite[Corollary 19.7]{gilmer}
(see also \cite{kang} and \cite{dobbs} for recent developments
in the lifting of chains to valuation rings),
and that a valuation domain is a Bezout domain and hence a
torsion free module over it is flat, see \cite[Theorem 63]{kaplansky}
and \cite[Ch.1, \S4, Proposition 3]{bourbakiac}.
\end{remarks}

\begin{corollary}
Let $A$ and $B$
be commutative rings and let $\varphi:A \ra B$ be a
ring homomorphism such that $I=\varphi^{-1} (IB)$ holds
for every ideal $I \subseteq A$.
Then $\dim \, B \geq \dim\, A$.
\end{corollary}

\proof
This is clear from the Theorem.
\qed

\begin{corollary}
Let $A$ be a commutative Noetherian
ring and let $B$ be an $A$-Algebra of finite type
such that every ideal of $A$ is contracted from $B$.
Then $g:\Spec \, B \ra \Spec \, A$ is submersive, i.e.
$\Spec\, A$ carries the quotient topology.
\end{corollary}

\proof
We have to show that a subset $W \subseteq \Spec \, A$ is open if
its preimage is open. Since $g$ is surjective,
we know that $W=g(g^{-1}(W))$, hence $W$ is constructible by
\cite[Th\'{e}or\`{e}me 7.1.4]{EGAI}.
For the openess it is therefore enough to show that it is closed
under generalization,
and this follows directly from our property:
let $\fop' \in W$ and let $\fop \subset \fop'$ be
a generalization. Let $\foq \subset \foq'$ be prime ideals lying
over them.
Then $\foq' \in g^{-1}(W)$ and since $g^{-1}(W) $ is open it is
closed under generalization, hence $\foq \in g^{-1}(W)$, and
this means $ \fop \in W$.
\qed

\medskip
It is easy to give an example of a direct summand
such that $\Spec\, B \ra \Spec \, A$
fulfills neither the going down nor the going up property.

\begin{example}
\label{example}
Let $K$ be a field and let the polynomial ring
$B=K[X,Y,Z]$ be $\ZZ$-graded
by $\deg \, X= \deg \, Y=1$,
$\deg \, Z =-1$.
Then the ring of degree zero is
$$A=B_0=K[XZ,YZ] \cong K[U,V]  \, .$$
$A$ is a direct summand in $B$, hence
the chain lifting property holds.

We consider the chain $(XZ) \subset (XZ,YZ)$ in $A$.
The principal prime ideal $ZB$ maps to $(XZ,YZ)$,
but no prime ideal $\subset ZB$ maps to $(XZ)$,
hence going down does not hold.

The prime ideal $(X,Y^2Z-1)B$ maps to $(XY)$.
But a prime ideal lying over $(XZ,YZ)$ must contain
either $ZB$ or $(X,Y)B$, hence also going up fails to hold.
\end{example}

\begin{remark}
A surjective (even bijective) mapping
between affine varieties may not
fulfill the chain lifting property, since there exist bijective
mappings which are not homeomorphisms.
\end{remark}


\end{document}